\ifx\shlhetal\undefinedcontrolsequence\let\shlhetal\relax\fi

\documentclass[11pt]{amsart}

\usepackage{amsmath}
\usepackage{amssymb}

\newtheorem{theorem}{Theorem}[section]
\newtheorem{claim}[theorem]{Claim}

\newtheorem{lemma}[theorem]{Lemma}

\newtheorem{corollary}[theorem]{Corollary}
\newtheorem{obs}[theorem]{Observation}

\theoremstyle{definition}
\newtheorem{definition}[theorem]{Definition}

\theoremstyle{remark}
\newtheorem{remark}[theorem]{Remark}

\newcount\skewfactor
\def\mathunderaccent#1#2 {\let\theaccent#1\skewfactor#2
\mathpalette\putaccentunder}
\def\putaccentunder#1#2{\oalign{$#1#2$\crcr\hidewidth
\vbox to.2ex{\hbox{$#1\skew\skewfactor\theaccent{}$}\vss}\hidewidth}}
\def\name{\mathunderaccent\tilde-3 }

\def\smallbox#1{\leavevmode\thinspace\hbox{\vrule\vtop{\vbox
   {\hrule\kern1pt\hbox{\vphantom{\tt/}\thinspace{\tt#1}\thinspace}}
   \kern1pt\hrule}\vrule}\thinspace}

\newcommand{\cf}{{\rm cf}}

\newcommand{\then}{{\underline{then}}}
\newcommand{\Then}{{\underline{Then}}}

\def\qedref#1{$\qed_{\reforiginal{#1}}$}

\title{A strong polarized relation}
\author{Shimon Garti}
\address{Institute of Mathematics
 The Hebrew University of Jerusalem
 Jerusalem 91904, Israel}
\email{shimon.garty@mail.huji.ac.il}

\author{Saharon Shelah}
\address{Institute of Mathematics
 The Hebrew University of Jerusalem
 Jerusalem 91904, Israel
 and  Department of Mathematics
 Rutgers University
 New Brunswick, NJ 08854, USA}
\email{shelah@math.huji.ac.il}
\urladdr{http://www.math.rutgers.edu/\char`\~shelah}

\thanks{First typed: December 2008 \newline Research supported by the United States-Israel Binational Science Foundation. This work is a part of the doctoral thesis of the first author, and publication 949 of the second author}

\begin{document}
\let\labeloriginal\label
\let\reforiginal\ref

\begin{abstract}
We prove that the strong polarized relation $\binom{\mu^+}{\mu} \rightarrow \binom{\mu^+}{\mu}^{1,1}_2$ is consistent with ZFC, for a singular $\mu$ which is a limit of measurable cardinals.
\end{abstract}

\maketitle

\newpage

\section{introduction}

The polarized relation $\binom{\alpha}{\beta} \rightarrow \binom{\gamma_0 \quad \gamma_1}{\delta_0 \quad \delta_1}^{1,1}$ asserts that for every coloring $c : \alpha \times \beta \rightarrow 2$ there are $A \subseteq \alpha$ and $B \subseteq \beta$ such that either ${\rm otp} (A) = \gamma_0, {\rm otp} (B) = \delta_0$ and $c \upharpoonright (A \times B) = \{0\}$ or ${\rm otp} (A) = \gamma_1, {\rm otp} (B) = \delta_1$ and $c \upharpoonright (A \times B) = \{1\}$. This relation was first introduced in \cite{MR0081864}, and investigated further in \cite{MR0202613}.

If $(\gamma_0, \delta_0) \neq (\gamma_1, \delta_1)$ then we get the so-called \emph{unbalanced form} of the relation. The \emph{balanced form} is the case $(\gamma_0, \delta_0) = (\gamma_1, \delta_1)$, and in this case we can write also $\binom{\alpha}{\beta} \rightarrow \binom{\gamma}{\delta}^{1,1}_2$ (stipulating $\gamma = \gamma_0 = \gamma_1$ and $\delta = \delta_0 = \delta_1$). With this shorthand, the notation $\binom{\alpha}{\beta} \rightarrow \binom{\gamma}{\delta}^{1,1}_\theta$ means the same thing, but the number of colors is $\theta$ instead of $2$.

From some trivialities and simple limitations, it follows that the case $\alpha = \mu^+$ and $\beta = \mu$ is interesting, for an infinite cardinal $\mu$. It is reasonable to distinguish between three cases - $\mu$ is a successor cardinal, $\mu$ is a limit regular cardinal (so it is a large cardinal) and $\mu$ is a singular cardinal; we concentrate in the latter case.

By a result of \v Cudnovski\u i  in \cite{MR0371655}, if $\mu$ is measurable then the relation $\binom{\mu^+}{\mu} \rightarrow \binom{\mu^+ \quad \alpha}{\mu \quad \mu}^{1,1}$ holds in ZFC for every $\alpha < \mu^+$ (see also \cite{MR1968607}, for discussion on weakly compact cardinals). In a sense, this is the best possible result, since we know that the assertion $\binom{\mu^+}{\mu} \nrightarrow \binom{\mu^+}{\mu}^{1,1}_2$ is valid under the GCH for every infinite cardinal $\mu$ (see \cite{williams}). This limitation gives rise to the following problem: Can one prove that the strong relation $\binom{\mu^+}{\mu} \rightarrow \binom{\mu^+}{\mu}^{1,1}_2$ is consistent with ZFC? For $\mu=\aleph_0$ the answer is yes. The same result holds for every supercompact cardinal $\mu$ (as we shall prove in a later work). But what happens if $\mu$ is singular?

We give here a positive answer. For a singular $\mu$ which is a limit of measurables, we can show that under some cardinal arithmetic assumptions (including the violation of the GCH, of course) one can get $\binom{\mu^+}{\mu} \rightarrow \binom{\mu^+}{\mu}^{1,1}_\theta$ for every $\theta < \cf (\mu)$. This result is stronger, on the one hand, than the balanced result $\binom{\mu^+}{\mu} \rightarrow \binom{\alpha}{\mu}^{1,1}_\theta$ which is proved in \cite{MR1606515} for every $\alpha < \mu^+$. On the other hand, the result there is proved in ZFC, whence the strong relations in this paper can not be proved in ZFC.

One can view this result as the parallel to the ordinary partition relation with respect to weakly compact cardinals. Recall that $\lambda\rightarrow(\theta,\kappa)^2$ means that for every coloring $c:[\lambda]^2\rightarrow 2$ there exists either $A\in[\lambda]^\theta$ so that $c\upharpoonright[A]^2=\{0\}$ or $B\in[\lambda]^\kappa$ such that $c\upharpoonright[B]^2=\{1\}$. 
We know that if $\lambda$ is inaccessible then $\lambda \rightarrow (\lambda, \alpha)^2$ for every $\alpha < \lambda$, but the strong (and balanced) relation $\lambda \rightarrow (\lambda)^2_2$ kicks $\lambda$ up in the chart of large cardinals, making it weakly compact. The result here is similar, replacing the ordinary partition relation by the polarized one.

Our notation is standard. We use the letters $\theta, \kappa, \lambda, \mu$ for infinite cardinals, and $\alpha, \beta, \gamma, \delta, \varepsilon, \zeta, i, j$ for ordinals. For a regular cardinal $\kappa$ we denote the ideal of bounded subsets of $\kappa$ by $J^{\rm bd}_\kappa$. For $A,B \subseteq \kappa$ we say $A \subseteq^* B$ when $A \setminus B$ is bounded in $\kappa$; the common usage of this symbol is for $\kappa = \aleph_0$, but here we apply it to uncountable cardinals.

Suppose $J$ is an ideal on $\kappa$.
The product $\prod \limits_{\varepsilon < \kappa} \lambda_\varepsilon / J$ is $\theta$-directed if every subset of cardinality \emph{less} than $\theta$ has an upper bound in the product (with respect to $<_J$). This applies also to products of partially ordered sets. For more information about cardinal arithmetic, the reader may consult \cite{MR1318912}.

For a measurable cardinal $\kappa$ and a normal ultrafilter $U$ on $\kappa$, let $\mathbb{Q}_U$ be the usual Prikry forcing. If $p=(t_p,A_p)\in \mathbb{Q}_U$ then $A_p$ is the pure component of $p$ and $t_p$ is the impure component. For infinite cardinals $\kappa, \lambda$ so that $\kappa < \lambda$ we denote by ${\rm Levy}(\kappa, \lambda)$ the Levy collapse of $\lambda$ to $\kappa$. This forcing notion consists of the partial functions $f : \kappa \rightarrow \lambda$ such that $|{\rm Dom}(f)| < \kappa$, ordered by inclusion. It collapses $\lambda$ to $\kappa$, and in general does not do any essential harm and does not change important things out of the interval $[\kappa, \lambda]$.

We adopt the convention that $p \leq q$ means $q$ gives more information than $p$ in forcing notions. We use the symbol $p \parallel_{\mathbb{P}} q$ in the sense that the conditions $p$ and $q$ are compatible in $\mathbb{P}$.

Throughout the paper, $\jmath_D : {\rm \bf{V}} \rightarrow M$ is the canonical elementary embedding of the universe into the transitive collapse $M$ of ${\rm \bf {V}}^\mu / D$ (where $D$ is a nonprincipal $\mu$-complete ultrafilter on $\mu$). $\mu$ is the critical point of $\jmath_D$, which means that $\mu$ is the first ordinal moved by $\jmath_D$. We shall use $\jmath$ instead of $\jmath_D$, when no confusion arises. The picture is as follows:

$$
\jmath : {\rm \bf{V}} \hookrightarrow {\rm \bf {V}}^\mu / D \cong M
$$

\par \noindent and we can treat $\jmath$ as a function from ${\rm \bf{V}}$ into $M$. We shall use the following basic result of Solovay, which asserts that if $\lambda$ is supercompact, $\tau \geq \lambda$ and $U_\tau$ is a fine and normal ultrafilter on $[\tau]^{<\lambda}$, then $\jmath_\tau (\lambda) > \tau$. We use a supercompact cardinal, but probably hyper-measurable cardinal suffices (as in \cite{MR1632081}, for example). We indicate, further, that being a limit of measurables (as we assume for our singular cardinal) can be weakend, and we hope to shed light on this subject in a subsequent work.

The paper is arranged in three sections. In the first one we prove the main result, in the second we deal with forcing preliminaries, and in the last one we deal with cardinal arithmetic theorems. \newline 

We thank the referees for their excellent work, the careful reading, corrections, clarifications and improvements.

\newpage 

\section{The combinatorial theorem}

\par \noindent We state the main result of the paper:

\begin{theorem}
\label{mt}
The main result. \newline 
Let $\mu$ be a singular cardinal, $\kappa = \cf (\mu)$ and $\theta < \kappa$. \newline 
Assume $2^\kappa < \cf(\lambda) \leq \lambda < \cf(\Upsilon)\leq \Upsilon \leq 2^\mu$.

Suppose $\mu$ is a limit of measurable cardinals, $\mu<\cf(\lambda)$, $\bar{\lambda} = \langle \lambda_\varepsilon : \varepsilon < \kappa \rangle$ is a sequence of measurables with limit $\mu$ so that $\kappa < \lambda_0$, $\prod \limits_{\varepsilon < \kappa} \lambda_\varepsilon / J^{\rm bd}_\kappa$ and $\prod \limits_{\varepsilon < \kappa} \lambda^+_\varepsilon / J^{\rm bd}_\kappa$ are $\cf(\Upsilon)$-directed, and $2^{\lambda_\varepsilon} = \lambda_\varepsilon^+$ for every $\varepsilon < \kappa$.

For every $\varepsilon < \kappa$ let $D_\varepsilon$ be a normal uniform ultrafilter on $\lambda_\varepsilon$, so the product $\prod \limits_{\varepsilon < \kappa} (D_\varepsilon, \subseteq^*) / J^{\rm bd}_\kappa$ is $\cf(\Upsilon)$-directed.

\Then\ the strong relation $\binom{\lambda}{\mu} \rightarrow \binom{\lambda}{\mu}^{1,1}_\theta$ holds.
\end{theorem}

\par \noindent \emph{Proof}. \newline 
We start with $\bar{\lambda}, \bar{D}$ as in the assumptions of the theorem, ensured by Claim \ref{directed0} and Theorem \ref{directed1} below. Given a coloring $c : \lambda \times \mu \rightarrow \theta$, we have to find a single color $i_* < \theta$ and two sets $A \in [\lambda]^\lambda, B \in [\mu]^\mu$ such that $c \upharpoonright (A \times B) = \{i_*\}$.

For every $\alpha < \lambda$ we would like to define the sequence of colors $\bar{i}_\alpha = \langle i_{\alpha, \varepsilon} : \varepsilon < \kappa \rangle$, so $i_{\alpha, \varepsilon} < \theta$ for every $\varepsilon < \kappa$. Suppose $\alpha < \lambda$ and $\varepsilon < \kappa$ are fixed. Since $D_\varepsilon$ is an ultrafilter, moreover, $D_\varepsilon$ is $\theta^+$-complete, there is an ordinal $i_{\alpha, \varepsilon} < \theta$ so that:

$$
A_{\alpha, \varepsilon} =_{\rm def} \{ \gamma < \lambda_\varepsilon : c(\alpha, \gamma) = i_{\alpha, \varepsilon} \} \in D_\varepsilon
$$

\par \noindent Let $\bar{A}_\alpha$ be the sequence $\langle A_{\alpha, \varepsilon} : \varepsilon < \kappa \rangle$, for every $\alpha < \lambda$. Without loss of generality, $A_{\alpha, \varepsilon} \cap \bigcup \limits_{\zeta < \varepsilon} \lambda_\zeta = \emptyset$ for every $\alpha < \lambda$ and $\varepsilon < \kappa$ (we can cut any initial segment of $A_{\alpha, \varepsilon}$ and still remain in the ultrafilter).

Recall that $\mu$ is limit of measurable cardinals, so in particular it is strong limit. Consequently, $\theta^\kappa < \mu < \lambda$, and even if $\lambda$ is singular we have $\theta^\kappa = 2^\kappa < \cf(\lambda) \leq \lambda$, so we have less than $\cf(\lambda)$ color-sequences of the form $\bar{i}_\alpha$ and we can choose $S_0 \subseteq \lambda$, $|S_0| = \lambda$ and a constant sequence $\langle i_\varepsilon : \varepsilon < \kappa \rangle$ so that:

$$
\alpha \in S_0 \Rightarrow \bar{i}_\alpha \equiv \langle i_\varepsilon : \varepsilon < \kappa \rangle
$$

\par \noindent Moreover, since $\theta < \kappa = \cf (\kappa)$ we can pick up an ordinal $i_* < \theta$ and a set $u \in [\kappa]^\kappa$ such that $\varepsilon \in u \Rightarrow i_\varepsilon \equiv i_*$. Without loss of generality, $u = \kappa$ (one may replace $\bar{\lambda}, \bar{D}$ by the sequences $\langle \lambda_\varepsilon : \varepsilon \in u \rangle$ and $\langle D_\varepsilon : \varepsilon \in u \rangle$).

The crucial step is the following: we choose a sequence of sets $\bar{A_*} = \langle A_\varepsilon^* : \varepsilon < \kappa \rangle$, $\bar{A_*} \in \prod \limits_{\varepsilon < \kappa} D_\varepsilon$, such that $A_\varepsilon^* \setminus A_{\alpha, \varepsilon}$ is bounded for (many and without loss of generality) each $\alpha \in S_0$ and every $\varepsilon < \kappa$. How can we ensure that such a sequence does exist? Well, each $\bar{A}_\alpha$ is a member in the product $\prod \limits_{\varepsilon < \kappa} (D_\varepsilon, \subseteq^*) / J^{\rm bd}_\kappa$. Since $|S_0| = \lambda < \Upsilon$ and by the $\Upsilon$-directness of the product, we can choose $\bar{A}_*$ such that:

$$
\alpha \in S_0 \Rightarrow \bar{A}_* \leq_{J^{\rm bd}_\kappa} \bar{A}_\alpha
$$

\par \noindent The meaning of the former is that $A_\varepsilon^* \setminus A_{\alpha, \varepsilon}$ is bounded for each $\alpha \in S_0$ (recall that the order of the product is reverse $\subseteq^*$). More precisely, $A_\varepsilon^* \setminus A_{\alpha, \varepsilon}$ is bounded for all large $\varepsilon$, but since $\kappa<\cf(\lambda)$ we can shrink $S_0$ and confine ourselves to a tail end of $\lambda_\varepsilon$-s.
We employ a similar argument to show that (after some shrinking of the set $S_0$) ${\rm sup} (A_\varepsilon^* \setminus A_{\alpha, \varepsilon})$ does not depend on $\alpha$. For this, define $g_\alpha \in \prod \limits_{\varepsilon < \kappa} \lambda_\varepsilon$ by $g_\alpha (\varepsilon) = {\rm sup} (A_\varepsilon^* \setminus A_{\alpha, \varepsilon}) < \lambda_\varepsilon$. We choose, in this way, just $\lambda$ functions. Since the product is $\Upsilon$-directed, there is $g_* \in \prod \limits_{\varepsilon < \kappa} \lambda_\varepsilon$ such that:

$$
\alpha \in S_0 \Rightarrow g_\alpha <_{J^{\rm bd}_\kappa} g_*
$$

\par \noindent Now define $j_\alpha = {\rm sup} \{ \varepsilon < \kappa : g_\alpha (\varepsilon) \geq g_* (\varepsilon) \} < \kappa$, for every $\alpha \in S_0$. $\kappa < \cf(\lambda) \leq \lambda$, so one can choose $S_1 \subseteq S_0$, $|S_1| = \lambda$ and an ordinal $j(*) < \kappa$ so that:

$$
\alpha \in S_1 \Rightarrow j_\alpha = j(*)
$$

\par \noindent Without loss of generality, $A^*_\varepsilon \cap [\bigcup \limits_{\zeta < \varepsilon} \lambda_\zeta, g_*(\varepsilon)) = \emptyset$, so one can verify that $(\forall \alpha \in S_1)(\forall \varepsilon \in [j(*), \kappa))(A^*_\varepsilon \subseteq A_{\alpha, \varepsilon})$, and we can construct now the desired sets $A$ and $B$. Define $A = S_1$, and $B = \bigcup \{ A^*_\varepsilon : j(*) \leq \varepsilon < \kappa \}$. Clearly, $A \in [\lambda]^\lambda$, and $B \in [\mu]^\mu$.

Suppose $\alpha \in A$ and $\beta \in B$. By the nature of $B$, there exists an ordinal $\varepsilon \in [j(*),\kappa)$ such that $\beta \in A^*_\varepsilon$, and since $A^*_\varepsilon \subseteq A_{\alpha, \varepsilon}$ we have $\beta \in A_{\alpha, \varepsilon}$, so $c(\alpha, \beta) = i_{\alpha, \varepsilon} = i_*$, and the relation $\binom{\lambda}{\mu} \rightarrow \binom{\lambda}{\mu}^{1,1}_\theta$ is established.

\hfill \qedref {mt}

\begin{corollary}
\label{mcoroll}
The strong polarized relation $\binom{\mu^+}{\mu} \rightarrow \binom{\mu^+}{\mu}^{1,1}_2$ is consistent with ZFC for some singular cardinal $\mu$.
\end{corollary}

\par\noindent\emph{Proof}.\newline 
As in Claim \ref{directed0} we prove the consistency of the conditions of Theorem \ref{mt} with $\kappa=\omega=\cf(\mu), \lambda=\mu^+$ and $\Upsilon=\mu^{++}=2^\mu$.

\hfill \qedref{mcoroll}

\begin{remark}
\label{st}
Denote by $\binom{\mu^+}{\mu} \rightarrow_{\rm st} \binom{\mu^+}{\mu}^{1,1}_\theta$ the assertion that for every coloring $c : \mu^+ \times \mu \rightarrow \theta$ there are $A$ and $B$ such that $A$ is a stationary subset of $\mu^+$, $B \in [\mu]^\mu$ and $c$ is constant on the cartesian product $A \times B$. Actually, our proof gives this relation.
\end{remark}

\newpage

\section{forcing preliminaries}

\par \noindent We need some preliminaries, before proving the main claim of the next section. First of all, we shall use a variant of Laver's indestructibility (see \cite{MR0472529}), making sure that a supercompact cardinal $\lambda$ will remain supercompact upon forcing with some prescribed properties. Let us start with the following definition:

\begin{definition}
\label{strategic}
Strategical completeness. \newline 
Let $\mathbb{P}$ be a forcing notion, $p \in \mathbb{P}$, and let $\mu$ be an infinite cardinal.
\begin{enumerate}
\item [$(a)$] The game $\Game_\mu(p, \mathbb{P})$ is played between two players, `com' and `inc'. It lasts $\mu$ moves. In the $\alpha$-th move, `com' tries to choose $p_\alpha \in \mathbb{P}$ such that $p \leq_{\mathbb{P}} p_\alpha$ and $\beta < \alpha \Rightarrow q_\beta \leq_{\mathbb{P}} p_\alpha$. After that, `inc' tries to choose $q_\alpha \in \mathbb{P}$ such that $p_\alpha \leq_{\mathbb{P}} q_\alpha$.
\item [$(b)$] `com' wins a play if he has a legal move for every $\alpha < \mu$.
\item [$(c)$] $\mathbb{P}$ is $\mu$-strategically complete if the player `com' has a winning strategy in the game $\Game_\mu(p, \mathbb{P})$ for every $p \in \mathbb{P}$.
\end{enumerate}
\end{definition}

\begin{claim}
\label{laver}
Indestructible supercompact and strategically completeness. \newline 
Let $\lambda$ be a supercompact cardinal in the ground model.
There is a forcing notion $\mathbb{Q}$ which makes $\lambda$ indestructible under every forcing $\mathbb{P}$ with the following properties:
\begin{enumerate}
\item [$(a)$] $\mathbb{P}$ is $\mu$-strategically complete for every $\mu < \lambda$,
\item [$(b)$] $\chi \geq \lambda$, and $\mathbb{P} \in \mathcal{H}(\chi)$,
\item [$(c)$] for some $\jmath : {\rm \bf{V}} \rightarrow M$ such that $\lambda = {\rm crit}(\jmath)$, $M^\chi \subseteq M$ and for every $G \subseteq \mathbb{P}$ which is generic over {\rm \bf{V}}, we have $M[G] \models "\{ \jmath(p) : p \in G \}$ has an upper bound in $\jmath(\mathbb{P})"$.
\end{enumerate}
\end{claim}

\par \noindent \emph{Proof}. \newline 
Basically, the proof walks along the line of \cite{MR0472529}, using Laver's diamond. In the crux of the matter, when Laver needs the $\lambda$-completeness, we employ requirement (c) above.

\hfill \qedref{laver}

\par \noindent We define now the `single step' forcing notion $\mathbb{Q}_{\bar{\theta}}$ to be used in the proof of Claim \ref{directed0}. This is called the $\bar{\theta}$-dominating forcing (it appears also in \cite{945}). We will use an iteration which consists, essentially, of these forcing notions:

\begin{definition}
\label{dominating}
The $\bar{\theta}$-dominating forcing. \newline 
Let $\lambda$ be a supercompact cardinal.
Suppose $\bar{\theta} = \langle \theta_\alpha : \alpha < \lambda \rangle$ is an increasing sequence of regular cardinals so that $2^{|\alpha|+\aleph_0} < \theta_\alpha < \lambda$ for every $\alpha < \lambda$. 
\begin{enumerate}
\item [$(\aleph)$] $p \in \mathbb{Q}_{\bar{\theta}}$ iff:
\begin{enumerate}
\item $p = (\eta, f) = (\eta^p, f^p)$,
\item $\ell g(\eta) < \lambda$,
\item $\eta \in \prod \{ \theta_\zeta : \zeta < \ell g(\eta) \}$,
\item $f \in \prod \{ \theta_\zeta : \zeta < \lambda \}$,
\item $\eta \triangleleft f$ (i.e., $\eta(\zeta)=f(\zeta)$ for every $\zeta<\ell g(\eta)$).
\end{enumerate}
\item [$(\beth)$] $p \leq_{\mathbb{Q}_{\bar{\theta}}} q$ iff ($p,q \in \mathbb{Q}_{\bar{\theta}}$ and) 
\begin{enumerate}
\item $\eta^p \trianglelefteq \eta^q$,
\item $f^p(\varepsilon) \leq f^q(\varepsilon)$, for every $\varepsilon < \lambda$.
\end{enumerate}
\end{enumerate}
\end{definition}

Notice that if $\ell g(\eta^p) \leq \varepsilon < \ell g(\eta^q)$ then $f^p(\varepsilon) \leq \eta^q(\varepsilon)$, since $f^p(\varepsilon) \leq f^q(\varepsilon) = \eta^q(\varepsilon)$. The purpose of $\mathbb{Q}_{\bar{\theta}}$ is to add (via the generic object) a dominating function in the product of the $\theta_\alpha$-s.

\begin{obs}
\label{cc}
Basic properties of $\mathbb{Q}_{\bar{\theta}}$. \newline 
Let $\mathbb{Q}_{\bar{\theta}}$ be the $\bar{\theta}$-dominating forcing (for the supercompact cardinal $\lambda$).
\begin{enumerate}
\item [$(a)$] $\mathbb{Q}_{\bar{\theta}}$ satisfies the $\lambda^+$-cc.
\item [$(b)$] $\mathbb{Q}_{\bar{\theta}}$ is $\mu$-strategically complete for every $\mu < \lambda$.
\end{enumerate}
\end{obs}

\par \noindent \emph{Proof}.
\begin{enumerate}
\item [$(a)$] If $p = (\eta, f^p), q = (\eta, f^q)$, define $f(\varepsilon) = {\rm max} \{ f^p(\varepsilon), f^q(\varepsilon) \}$ for every $\varepsilon < \lambda$, and then $r = (\eta, f)$. Clearly, $r \in \mathbb{Q}_{\bar{\theta}}$ and $p,q \leq r$. So the cardinality of an antichain does not exceed the number of possible $\eta$-s, which is $\lambda$ since $\ell g(\eta) < \lambda$ and $\lambda^{< \lambda} = \lambda$.
\item [$(b)$] Assume $\mu < \lambda$ and $p \in \mathbb{Q}_{\bar{\theta}}$. Let $\Game_\mu (p, \mathbb{Q}_{\bar{\theta}})$ be the game defined in \ref{strategic}. We shall find a winning strategy for `com'. In the first stage, `com' may choose $p_0 = p$, and from now on `com' needs to deal only with the $q_\gamma$-s (and being above $p$ follows, since $p \leq_{\mathbb{Q}_{\bar{\theta}}} q_\gamma$ for each $\gamma$). So assume $\beta < \mu$ and $q_\gamma$ was already chosen (by `inc') for $\gamma < \beta$. Notice that $\langle q_\gamma : \gamma < \beta \rangle$ is an increasing sequence of conditions in $\mathbb{Q}_{\bar{\theta}}$. Define:

$$
\eta^{p_\beta} = \bigcup \{ \eta^{q_\gamma} : \gamma < \beta \}
$$

\par \noindent Since $\beta < \lambda = \cf(\lambda)$ and $\ell g(\eta^{q_\gamma}) < \lambda$ for every $\gamma < \beta$, we know that $\ell g(\eta^{p_\beta}) < \lambda$. Note that $\eta^{q_\gamma} \trianglelefteq \eta^{p_\beta}$ for every $\gamma < \beta$. \newline 
Now, for $\varepsilon < \ell g(\eta^{p_\beta})$ set $f^{p_\beta}(\varepsilon) = \eta^{p_\beta}(\varepsilon)$, and for $\ell g(\eta^{p_\beta}) \leq \varepsilon < \lambda$ set $f^{p_\beta}(\varepsilon) = {\rm sup} \{ f^{q_\gamma}(\varepsilon) : \gamma < \beta \}$. We may assume, without loss of generality, that $\beta<\ell g(\eta^{p_\beta})$ (if not, set $f^{p_\beta}(\varepsilon)=0$ for every $\varepsilon\in[\ell g(\eta^{p_\beta}),\beta]$).
Hence $f^{p_\beta}(\varepsilon)$ is well defined, since $\alpha < \theta_\alpha$ for every $\alpha < \lambda$. Notice also that $\eta^{p_\beta} \triangleleft f^{p_\beta}$. Finally, set $p_\beta = (\eta^{p_\beta}, f^{p_\beta})$. Clearly, $p_\beta \in \mathbb{Q}_{\bar{\theta}}$ and $q_\gamma \leq_{\mathbb{Q}_{\bar{\theta}}} p_\beta$ for every $\gamma < \beta$, so we are done.
\end{enumerate}

\hfill \qedref{cc}

\begin{remark}
\label{comp}
Despite the fact that $\mathbb{Q}_{\bar{\theta}}$ is $\mu$-strategically complete for every $\mu<\lambda$, it is not $\lambda$-complete. We indicate that for every $\mu < \lambda$ there is a dense subset which is $\mu$-complete, but no dense subset which is $\mu$-complete simultaneously for every $\mu < \lambda$. So we will have to employ claim \ref{laver} instead of the original theorem of Laver.
\end{remark}

\hfill \qedref{comp}

\par \noindent Having the basic component, we would like to iterate the $\bar{\theta}$-dominating forcing. We shall use a $(< \lambda)$-support, aiming to take care of all the increasing sequences of the form $\bar{\theta}$ with limit $\lambda$. We need the following:

\begin{definition}
\label{iteration}
The iteration. \newline 
Let $\lambda$ be a supercompact cardinal, and $\lambda<\cf(\Upsilon)\leq\Upsilon$.
Let $\mathbb{P}_\Upsilon$ be the $(< \lambda)$-support iteration $\langle \mathbb{P}_\alpha, \name{\mathbb{Q}}_\beta : \alpha \leq \Upsilon, \beta < \Upsilon \rangle$, where each $\name{\mathbb{Q}}_\beta$ is (a $\mathbb{P}_\beta$-name of) the forcing $\mathbb{Q}_{\bar{\theta}}$ with respect to some $\bar{\theta}$ as in definition \ref{dominating}, so that each $\bar{\theta}$ appears at some stage of the iteration.
\end{definition}

\par \noindent We would like to show that the nice properties of each component ensured by \ref{cc} are preserved in the iteration. Now, the strategical completeness is preserved, but the chain condition may fail. Nevertheless, in the case of the dominating forcing $\mathbb{Q}_{\bar{\theta}}$ it holds:

\begin{definition}
\label{strongcc}
Linked forcing notions. \newline 
Let $\mathbb{P}$ be a forcing notion.
$\mathbb{P}$ is $\lambda$-2-linked when for every subset of conditions $\{ p_\alpha : \alpha < \lambda^+ \} \subseteq \mathbb{P}$ there are $C,h$ such that:
\begin{enumerate}
\item $C$ is a closed unbounded subset of $\lambda^+$
\item $h : \lambda^+ \rightarrow \lambda^+$ is a regressive function
\item for every $\alpha, \beta \in C$, if $\cf(\alpha) = \cf(\beta) = \lambda$ and $h(\alpha) = h(\beta)$ then $p_\alpha \parallel p_\beta$; moreover, $p_\alpha \cup p_\beta$ is a least upper bound
\end{enumerate}
\end{definition}

\begin{remark}
\label{strongness}
If $\mathbb{Q}$ is $\lambda$-2-linked, \then\ $\mathbb{Q}$ is $\lambda^+$-cc.
\end{remark}

\begin{lemma}
\label{preservation}
Preservation of the $\lambda$-2-linked property. \newline 
Assume $\lambda = \lambda^{< \lambda}$, $\mathbb{P}$ is a $(< \lambda)$-support iteration so that every component is $\chi$-strategically complete for every $\chi < \lambda$ and $\lambda$-2-linked. \newline 
\Then\ $\mathbb{P}$ is also $\lambda$-2-linked (and consequently $\lambda^+$-cc).
\end{lemma}

\par \noindent \emph{Proof}. \newline 
As in \cite{MR0505492}, with the minor changes for $\lambda$ instead of $\aleph_1$.

\hfill \qedref{preservation}

\begin{obs}
\label{dominatcc}
The dominating forcing $\mathbb{Q}_{\bar{\theta}}$ is $\lambda$-2-linked.
\end{obs}

\par \noindent \emph{Proof}. \newline 
Suppose $\{ p_\alpha : \alpha < \lambda^+ \} \subseteq \mathbb{Q}_{\bar{\theta}}$. Without loss of generality there exists $\eta$ so that $\eta^{p_\alpha} \equiv \eta$ for every $\alpha < \lambda^+$ (since $\ell g(\eta) < \lambda$ and $\lambda^{< \lambda} = \lambda < \lambda^+$). Now choose any club $C$ and regressive $h$, upon noticing that two conditions with the same stem are compatible.

\hfill \qedref{dominatcc}

\begin{lemma}
\label{prodcof}
The high cofinality. \newline 
Suppose $\lambda$ is supercompact, $\cf(\Upsilon) > \lambda$, $\bar{\theta}$ and $\mathbb{P} = \mathbb{P}_\Upsilon$ are the sequence of regular cardinals and iteration defined above. 
 
\Then\ $\cf(\prod \limits_{\alpha < \lambda} \theta_\alpha, <_{J^{\rm {bd}}_\lambda}) = \cf(\Upsilon)^{{\rm \bf V}^\mathbb{P}}$.
\end{lemma}

\par \noindent \emph{Proof}. \newline 
Let $\Upsilon = \lambda^{++}$ (the proof of the general case is just the same). 
For proving that $\cf(\prod \limits_{\alpha < \lambda} \theta_\alpha, <_{J^{\rm {bd}}_\lambda}) \leq \lambda^{++}$ we introduce a cofinal subset (in the product) of cardinality $\lambda^{++}$. Moreover, the cofinal subset will be a dominating one (i.e., each $\name{g}_\beta$ below dominates all the old functions), and consequently $\cf(\prod \limits_{\alpha < \lambda} \theta_\alpha, <_{J^{\rm {bd}}_\lambda}) = \lambda^{++}$.

For each $\beta < \lambda^{++}$ let $\name{G}_\beta \subseteq \name{\mathbb{Q}}_\beta$ be generic, and set $\name{g}_\beta = \bigcup \{ \eta^p : p \in \name{G}_\beta \}$. Now,  $\name{\mathbb{P}}_{\beta+1} \models ``\name{g}_\beta \in \prod \limits_{\alpha < \lambda} \theta_\alpha$ and $\name{g}_\beta$ is a dominating function". To see this, define the following set for every $g \in \prod \limits_{\alpha < \lambda} \theta_\alpha$:

$$
\mathcal{I}_g = \{ (\eta, f) \in \name{\mathbb{Q}}_\beta : \forall \varepsilon \in [\ell g(\eta), \lambda) \quad g(\varepsilon) \leq f(\varepsilon) \}
$$

\par \noindent One verifies that $\mathcal{I}_g$ is a dense open set for every $g \in \prod \limits_{\alpha < \lambda} \theta_\alpha$, so if $G$ is generic then $G \cap \mathcal{I}_g \neq \emptyset$ for every $g \in \prod \limits_{\alpha < \lambda} \theta_\alpha$. Consequently, $\Vdash_{\name{\mathbb{P}}_{\beta+1}} "g \leq_{J^{\rm {bd}}_\lambda} \name{g}_\beta"$. Take a look at $\{ \name{g}_\beta : \beta < \lambda^{++} \}$. We claim (working in ${\rm \bf{V}}^{\mathbb{P}}$) that this set is cofinal in the product.

For showing this, notice that $\langle \mathbb{P}_\alpha : \alpha \leq \lambda^{++} \rangle$ is $\lessdot$-increasing, so $\alpha < \beta < \lambda^{++} \Rightarrow {\rm \bf{V}}^{\mathbb{P}} \models "\name{g}_\alpha \leq_{J^{\rm {bd}}_\lambda} \name{g}_\beta"$, and by the nature of these objects we know that every function in the product is bounded by one of them.

\hfill \qedref{prodcof}

\begin{lemma}
\label{indestructible}
The property of being indestructible. \newline 
The iteration $\mathbb{P}$ satisfies demand $(c)$ in claim \ref{laver}.
\end{lemma}

\par \noindent \emph{Proof}. \newline 
Let $G \subseteq \mathbb{P}$ be generic over ${\rm \bf{V}}$. Let $\jmath$ be a $\chi$-supercompact elementary embedding of ${\rm \bf{V}}$ into $M$ with critical point $\lambda$, and let $\Upsilon = \lambda^{++}$. We may assume that $\chi \geq \Upsilon$. We define a condition $q \in \jmath(\mathbb{P})$, and we shall prove that $q$ is an upper bound (in the forcing notion $\jmath(\mathbb{P})$ which belongs to $M[G]$) for $\{ \jmath(p) : p \in G \}$.

Set ${\rm Dom}(q) = \{ \jmath(\alpha) : \alpha < \Upsilon \}$. By saying, below, that $\eta^p$ is an object we mean that it is not just a name.
For every $\alpha < \Upsilon$ let $q(\jmath(\alpha)) = (\eta^\alpha, \name{f}^\alpha)$ where:

$$
\eta^\alpha = \bigcup \{ \eta^{p(\alpha)} : p \in G, \alpha \in {\rm Dom}(p), \eta^p {\rm \ is \ an \ object} \}
$$

\par \noindent and for $\gamma \geq \ell g(\eta^\alpha)=\lambda$, set:

$$
f^\alpha(\gamma) = {\rm sup} \{ \jmath(f^{p(\alpha)})(\gamma) : p \in G, \alpha \in {\rm Dom}(p) \}
$$

\par \noindent Clearly, $q$ is an upper bound for $\{ \jmath(p) : p \in G \}$ in $\jmath(\mathbb{P})$, provided that $q$ is well defined. For this, notice that $|{\rm Dom}(q)| < \jmath(\lambda)$ since $\Upsilon < \jmath(\lambda)$, $^\chi M \subseteq M$, and ${\rm Dom}(q) = \{ \jmath(\alpha) : \alpha < \Upsilon \}$.

We also must show that $f^\alpha(\gamma)$ is well defined. Notice that for every $\bar{\theta}$ which proceeds fast enough (i.e., $\alpha < \lambda \Rightarrow 2^{|\alpha|+\aleph_0} < \theta_\alpha < \lambda$ as in \ref{dominating}) and for each $\gamma \geq \lambda$ we have $\jmath(\bar{\theta})_\gamma > \Upsilon$. Also, if $\alpha < \Upsilon$ then $M[G] \models |\{ f^{p(\alpha)} : p \in G_{\mathbb{P}} \}| < \jmath(\bar{\theta})_\gamma$ for $\gamma \geq \lambda$. Consequently, $f^\alpha(\gamma)$ is bounded in $\jmath(\bar{\theta})_\gamma$, hence well defined, so we are done.

\hfill \qedref{indestructible}

\newpage

\section{Cardinal arithmetic assumptions}

\par \noindent We phrase two theorems, which we shall prove in this section. The first one asserts that there exists (i.e., by forcing) a singular cardinal, limit of measurable cardinals, with some properties imposed on the product of these measurables and their successors. Related works, in this light, are \cite{MR1632081} and \cite{MR1245523}. The second theorem deals with properties of the product of normal (uniform) ultrafilters on these cardinals. \newline 
We start with the following known fact:

\begin{lemma}
\label{prikry}
Cofinality preservation under Prikry forcing. \newline 
Let $U$ be a normal (uniform) ultrafilter on a measurable cardinal $\mu$. \newline 
Let $\mathbb{Q}_U$ be the Prikry forcing (with respect to $\mu$ and $U$) and $\langle \vartheta_n : n \in \omega \rangle$ the Prikry sequence.

Suppose $\theta = \cf(\theta) \neq \mu$, $F : \mu \rightarrow {\rm Reg} \cap \mu$, $F(\alpha)>\alpha$ for every $\alpha<\mu$ and $\cf(\prod \limits_{\alpha < \mu} F(\alpha) / U) = \theta$ as exemplified by $\bar{g} = \langle g_\varepsilon / U : \varepsilon < \theta \rangle$ in ${\rm \bf {V}}$. Let $\bar{h}=\langle h_\varepsilon : \varepsilon < \theta \rangle$ be the restriction of $\bar{g}$ to the Prikry sequence, i.e., $h_\varepsilon(\vartheta_0)=0$ and $h_\varepsilon\upharpoonright\{\vartheta_{n+1}:n\in\omega\}= g_\varepsilon\upharpoonright\{\vartheta_{n+1}:n\in\omega\}$ for every $\varepsilon<\theta$.

\Then\ $\cf(\prod \limits_{n < \omega} F(\vartheta_n) / J_\omega^{\rm bd}) = \theta$, as exemplified by $\bar{h}$ in ${\rm \bf{V}}^{\mathbb{Q}_U}$.
\end{lemma}

\par \noindent \emph{Proof}. \newline 
Let $A$ be any member of $U$. We claim that $\vartheta_{i+1}\in A$ for almost every $i\in\omega$ (i.e., except a finite set). For this, let $D_A$ be $\{p\in\mathbb{Q}_U:A\supseteq A_p\}$ (recall that $A_p$ is the pure component of the condition $p$). Let $G$ be a generic subset of $\mathbb{Q}_U$.
Since $D_A$ is open and dense, one can pick a condition $p\in D_A\cap G$. Let $i_p$ be the maximal natural number so that $\vartheta_{i_p}\in t_p$. Consequently, $p\Vdash(\forall i\in[i_p,\omega))(\vartheta_{i_p}\in A)$.

Now suppose $p\Vdash (\name{f}:\omega\rightarrow\mu)\wedge (\bigwedge\limits_{i<\omega}\name{f}(i)<F(\name{\vartheta}_{i+1}))$. We claim that there exists a condition $q\geq p$ and a function $g\in{\rm \bf V}$ so that $g:\mu\rightarrow\mu, \bigwedge\limits_\lambda g(\lambda)<F(\lambda)$ and $q\Vdash \name{f}(i)<g(\lambda)$ whenever $\name{\vartheta}_{i+1}=\lambda$ (more precisely, this holds for every large enough $i$ since for every measure one set $A$, $\vartheta_{i+1}\in A$ for all large $i$). For this claim, let $A_p$ be the pure component of $p$. For each $\lambda\in A_p$ define:

$$
T_{p,\lambda}=\{t:\exists A' \in U, p\leq(t,A'), {\rm max}(t)=\lambda\}
$$

For every $t\in T_{p,\lambda}$ we choose a condition $q_{t,\lambda}$ so that:
\begin{enumerate}
\item [$(a)$] $q_{t,\lambda}$ is of the form $(t,A')$, so forces the value $\lambda$ to $\name{\vartheta}_{|t\cap\lambda|}$
\item [$(b)$] $q_{t,\lambda}$ forces a value to $\name{f}_{|t\cap\lambda|}$ which is an ordinal below $F(\lambda)$
\end{enumerate}

Denote this ordinal by $g_p(t,\lambda)$. Now we define a condition $q=(s,A)$ as follows. $s=t_p$, and $A=A_q$ is the following set:

$$
A=\{\lambda\in A_p:\forall\lambda_1\in\lambda\cap A_p, \forall t\in T_{p,\lambda}, \lambda\in A_{q_{t,\lambda_1}}\}
$$

We shall show that $q\in\mathbb{Q}_U$. $A\subseteq A_p$, hence ${\rm max}(t_p)<{\rm min}(A)$. $A\in U$ since $A$ is the diagonal intersection of $\mu$ members from $U$. To verify this, set $B_\lambda=\bigcap\{A_{q_{t,\lambda}}:t\in T_{p,\lambda}\}$, for every $\lambda\in A_p$. Now $B_\lambda\in U$ for every $\lambda\in A_p$, as an intersection of at most $\lambda=|[\lambda]^{<\omega}|$ members from $U$ (recall that $U$ is $\mu$-complete). Since $A=\Delta\{B_\lambda:\lambda\in A_p\}$ we know that $A\in U$, hence $q\in\mathbb{Q}_U$.

Clearly, $p\leq_{\mathbb{Q}_U}q$. Let us define $g:\mu\rightarrow\mu$ by $g(\lambda) = {\rm sup}\{g_p(t,\lambda)+1:t\in T_{p,\lambda}\}$ if $\lambda\in A_p$, and $g(\lambda)=0$ otherwise. Notice that $g(\lambda)<F(\lambda)$, since $\lambda<F(\lambda), F(\lambda)$ is regular, $g_p(t,\lambda)<F(\lambda)$ for every $t\in T_{p,\lambda}$ and $|T_{p,\lambda}|\leq\lambda$. 
It follows that $q\Vdash\name{f}(i)<g(\name{\vartheta}_{i+1})$ (for almost every $i$ hence without loss of generality for every $i$), as required.

We conclude that if $p\Vdash \name{f}\in\prod\limits_{i<\omega}F(\vartheta_i)$ then one can find a function $g\in\prod\limits_{\alpha<\mu}F(\alpha)$, an ordinal $j<\omega$ and a condition $q\geq p$ such that $q\Vdash 
\bigwedge\limits_{i\in[j,\omega)}\name{f}(i)<g(\name{\vartheta}_{i+1})$. Equipped with this property, we can accomplish the proof of the lemma.

For every $\varepsilon<\theta$ set $h_\varepsilon=g_\varepsilon\upharpoonright\{\vartheta_{i+1}:i<\omega\}\cup
\langle\vartheta_0,0\rangle$ and collect these functions to the sequence $\bar{h}= \langle h_\varepsilon:\varepsilon<\theta\rangle$. We claim that $\bar{h}$ is a cofinal sequence in the product $(\prod\limits_{n\in\omega}F(\vartheta_n),<_{J_\omega^{\rm bd}})$ (in ${\rm \bf{V}}^{\mathbb{Q}_U}$). Assume $p\Vdash \name{f}\in\prod\limits_{i<\omega}F(\vartheta_i)$, and let $g,q$ be as above. Pick an ordinal $\varepsilon<\theta$ so that $g<_U g_\varepsilon$. It means that $B_{\varepsilon,g}=\{\gamma<\mu:g(\gamma)<g_\varepsilon(\gamma)\}\in U$. By the beginning of the proof, $\vartheta_{i+1}\in B_{\varepsilon,g}$ for almost every $i$, hence $q\Vdash\name{f}<_{J_\omega^{\rm bd}}h_\varepsilon$, and we are done.

\hfill \qedref{prikry}

\begin{remark}
\label{mmmagidor}
The same proof works for Magidor's forcing, upon replacing $\omega$ by $\kappa=\cf(\mu)$. In the proof above we demanded $q\Vdash\name{f}(i)<g(\name{\vartheta}_{i+1})$ (and not $\name{\vartheta}_i$), so it works also for Magidor's forcing (in contrary to Prikry forcing, in Magidor's forcing we encounter limit points in the cofinal sequence).
\end{remark}

\par \noindent We can state now the main claim of this section:

\begin{claim}
\label{directed0}
The main claim. \newline 
Starting with a supercompact cardinal, one can force the existence of a singular cardinal $\mu > \cf(\mu) = \kappa$, limit of measurables $\bar{\lambda} = \langle \lambda_\varepsilon : \varepsilon < \kappa \rangle$, such that both $\prod \limits_{\varepsilon < \kappa} \lambda_\varepsilon / J^{\rm bd}_\kappa$ and $\prod \limits_{\varepsilon < \kappa} \lambda^+_\varepsilon / J^{\rm bd}_\kappa$ are $\cf(\Upsilon)$-directed (for some $\Upsilon \in [\mu^{++},2^\mu), \cf(\Upsilon)\geq\mu^{++}$), and $2^{\lambda_\varepsilon} = \lambda_\varepsilon^+$ for every $\varepsilon < \kappa$.
\end{claim}

\par \noindent \emph{Proof}. \newline 
We shall prove the claim for the specific case of $\kappa=\omega$. The arguments can be generalized upon using Magidor's forcing instead of $\mathbb{Q}_U$ below.
Let $\mu$ be a supercompact cardinal. Begin with the variant of Laver's forcing, ensured by \ref{laver} above. Let ${\rm Levy}(\mu^+, 2^\mu)$ follow Laver's forcing, so $2^\mu = \mu^+$ and $\mu$ remains supercompact (notice that ${\rm Levy}(\mu^+, 2^\mu)$ is $\mu$-directed-closed). Use $\mathbb{P}$ from definition \ref{iteration} to follow the composition of the preparatory Laver forcing and ${\rm Levy}(\mu^+, 2^\mu)$.

By observation \ref{cc} we know that $\mathbb{P}$ is $\chi$-strategically complete for every $\chi < \mu$ (recall that an iteration keeps this property, provided that each stage satisfies it) and also $\mu^+$-cc (by lemma \ref{preservation} and observation \ref{dominatcc} above). By claim \ref{laver} we know that $\mu$ is still supercompact after forcing with $\mathbb{P}$. It should be stretched that $\mathbb{P}$ forces $2^\mu>\mu$, as it creates a $\mu^{++}$-directed product.

Choose a sequence of measurable cardinals $\langle \lambda_\varepsilon : \varepsilon < \mu \rangle$, so that $\mu$ is the limit of the sequence. Notice that both $\langle \lambda_\varepsilon : \varepsilon < \mu \rangle$ and $\langle \lambda_\varepsilon^+ : \varepsilon < \mu \rangle$ fit the definition \ref{dominating} (hence appear at some stage of the iteration). Without loss of generality $\varepsilon < \mu \Rightarrow 2^{\lambda_\varepsilon} = \lambda_\varepsilon^+$ (recall ${\rm Levy}(\mu^+, 2^\mu)$ upon noticing that the local GCH on $\mu$ reflects down to enough measurables below, and the iteration $\mathbb{P}$ does not affect the measurability of the cardinals below $\mu$, since it does not add new bounded subsets).

Let $U$ be a normal (uniform) ultrafilter on $\mu$ in ${\rm \bf{V}}^{\mathbb{P}}$. Let $\mathbb{Q}_U$ be the Prikry forcing applied to $\mu$, adding the cofinal Prikry sequence $\langle \vartheta_n : n < \omega \rangle$. In ${\rm \bf{V}}^{\mathbb{P} \ast \mathbb{Q}_U}$ we know that $\cf(\mu) = \aleph_0$. We indicate that using Magidor's forcing (from \cite{MR0465868}), we can get a similar result for $\cf(\mu) = \kappa > \aleph_0$.

Now, if $\bar{\theta}$ is an increasing sequence as in definition \ref{dominating}, we know that $\cf(\prod \limits_{n < \omega} \theta_{\vartheta_n} / J^{{\rm bd}}_\omega)^{{\rm \bf{V}}^{\mathbb{P}\ast\mathbb{Q}_U}} = \cf(\prod \limits_{\alpha < \mu} \theta_\alpha / U)^{{\rm \bf{V}}^{\mathbb{P}}} = \cf(\prod \limits_{\alpha < \mu} \theta_\alpha / J^{{\rm bd}}_\mu)^{{\rm \bf{V}}^{\mathbb{P}}} = \cf(\Upsilon)$ (by \ref{prikry} and \ref{prodcof} above, and the second equality follows from the fact that we have here true cofinality, which is preserved under extending the ideal). Apply it to the sequences $\langle \lambda_{\vartheta_\varepsilon} : \varepsilon < \kappa \rangle$ and $\langle \lambda_{\vartheta_\varepsilon}^+ : \varepsilon < \kappa \rangle$ so the proof is complete. 

\hfill \qedref{directed0}

\begin{remark}
\label{more}
We have used the main theorem for proving a combinatorial result, but we indicate that it can serve for other problems as well (by describing an extreme situation of cardinal arithmetic).
\end{remark}

\begin{theorem}
\label{directed1}
Let $\mu > \cf (\mu) = \kappa$ be singular, limit of measurables. \newline Let $\bar{\lambda} = \langle \lambda_\varepsilon : \varepsilon < \kappa \rangle$ be a sequence of measurable cardinals, which tends to $\mu$. Assume $2^{\lambda_\varepsilon} = \lambda_\varepsilon^+$ for every $\varepsilon < \kappa$, and $\prod \limits_{\varepsilon < \kappa} \lambda^+_\varepsilon / J^{\rm bd}_\kappa$ is $\Upsilon$-directed (for some $\Upsilon \in [\mu^{++}, 2^\mu)$).

\Then\ for every sequence $\bar{D} = \langle D_\varepsilon : \varepsilon < \kappa \rangle$ such that $D_\varepsilon$ is a normal (hence $\lambda_\varepsilon$-complete) uniform ultrafilter on $\lambda_\varepsilon$ for every $\varepsilon < \kappa$, the product $\prod \limits_{\varepsilon < \kappa} (D_\varepsilon, \subseteq^*) / J^{\rm bd}_\kappa$ is $\Upsilon$-directed.
\end{theorem} 

\par \noindent \emph{Proof}. \newline 
For every $\varepsilon < \kappa$ choose a normal ultrafilter $D_\varepsilon$ on $\lambda_\varepsilon$ (recall that each $\lambda_\varepsilon$ is measurable). First we claim that $(D_\varepsilon, \subseteq^*_\varepsilon)$ is $\lambda_\varepsilon^+$-directed for every $\varepsilon < \kappa$. So we have to show that for every collection of less than $\lambda_\varepsilon^+$ sets from $D_\varepsilon$ we can find a set in the ultrafilter which is almost included in every member of the collection.

Indeed, if $\{ S_\beta : \beta < \delta \}$ is such a collection (and without loss of generality $\delta \leq \lambda_\varepsilon$), define $S = \Delta \{ S_\beta : \beta < \delta \}$, and by the normality of $D_\varepsilon$ we know that $S \in D_\varepsilon$. Since $S \subseteq^*_\varepsilon S_\beta$ for every $\beta < \delta$ (by the very definition of the diagonal intersection), we are done.

Second we claim that for every $\varepsilon < \kappa$ we can find a $\subseteq^*_\varepsilon$-decreasing sequence of sets from $D_\varepsilon$ of the form $\langle S_{\varepsilon,\alpha} : \alpha < \lambda_\varepsilon^+ \rangle$, such that for every $B \in D_\varepsilon$ there is $\alpha < \lambda_\varepsilon^+$ so that $S_{\varepsilon,\alpha} \subseteq^*_\varepsilon B$. This is justified by the assumption that $2^{\lambda_\varepsilon} = \lambda_\varepsilon^+$, so we can enumerate the members of $D_\varepsilon$ by $\{ A_\gamma : \gamma < \lambda^+_\varepsilon \}$. For every $\gamma < \lambda^+_\varepsilon$, choose an enumeration of the collection $\{ A_\beta : \beta < \gamma \}$ as $\{ A_\alpha : \alpha < \gamma' \}$ such that $\gamma' \leq \lambda_\varepsilon$ (the new enumeration is needed whenever $\gamma > \lambda_\varepsilon$, and we want to arrange our sets below $\lambda_\varepsilon$). Define $S_\gamma = \Delta \{ A_\alpha : \alpha < \gamma' \}$. Now, for every $B \in D_\varepsilon$ find an ordinal $\delta$ above the index of $B$ in the enumeration of the members of $D_\varepsilon$, and $S_{\varepsilon,\delta} \subseteq^*_\varepsilon B$ as required.

Now we can show that the product $\prod \limits_{\varepsilon < \kappa} (D_\varepsilon, \subseteq^*) / J^{\rm bd}_\kappa$ is $\Upsilon$-directed. Assume $A \subseteq \prod \limits_{\varepsilon < \kappa} (D_\varepsilon, \subseteq^*)$, and $|A| = \Upsilon' < \Upsilon$. A typical member $\bar{C}_\alpha \in A$ is a sequence of the form $\langle C_\varepsilon^\alpha : \varepsilon < \kappa \rangle$ such that $C_\varepsilon^\alpha \in D_\varepsilon$ for every $\varepsilon < \kappa$. Choose an enumeration $\{ \bar{C}_\alpha : \alpha < \Upsilon' \}$ of $A$.

For each $\alpha < \Upsilon'$ we assign a vector $\bar{j}_\alpha = \langle j^\alpha_\varepsilon : \varepsilon < \kappa \rangle$ in the product $\prod \limits_{\varepsilon < \kappa} \lambda_\varepsilon^+$ as follows. For $\varepsilon < \kappa$ let $j_\varepsilon^\alpha$ be the index of the set $C_\varepsilon^\alpha$ in the enumeration of the members of $D_\varepsilon$ mentioned above. Define, now, the following set:

$$
A' = \{ \bar{j}_\alpha : \alpha < \Upsilon' \} \subseteq \prod \limits_{\varepsilon < \kappa} \lambda^+_\varepsilon
$$

\par \noindent Clearly, $|A'| < \Upsilon$, and we assume that this product is $\Upsilon$-directed, so we can choose a member $\bar{j} \in \prod \limits_{\varepsilon < \kappa} \lambda_\varepsilon^+$ which is an upper bound of $A'$. It means that $\alpha < \Upsilon' \Rightarrow \bar{j}_\alpha \leq_{J^{\rm bd}_\kappa} \bar{j}$. \newline 
$\bar{j}$ produces a member $\bar{C}$ in the product $\prod \limits_{\varepsilon < \kappa} (D_\varepsilon, \subseteq^*)$, as follows: for each $\varepsilon < \kappa$ we define $C_\varepsilon = S_{\varepsilon,j_\varepsilon}$, and then $\bar{C} = \langle C_\varepsilon : \varepsilon < \kappa \rangle$. Now $\bar{C}$ is an upper bound for the set $A$, and the proof is complete.

\hfill \qedref{directed1}

The theorems above (and consequently, the strong relation that was proved in the first section) require the existence of a supercompact cardinal in the ground model. The combinatorial result is proved on a singular cardinal, limit of measurables. It is plausible to get similar results below the first measurable cardinal, and with weaker assumption than supercompact. We indicate that some forcing which kills the measurability but keeps enough properties of the product $\prod \limits_{\varepsilon < \kappa} (D_\varepsilon, \subseteq^*) / J_\kappa^{\rm bd}$ is required. We hope to continue this subject in a subsequent paper.

\newpage 
\bibliographystyle{amsplain}
\bibliography{arlist}

\end{document}